\input amstex
\input epsf

 \documentstyle{amsppt}
 \magnification=1200
 \vsize=8.5truein
 \hsize= 6.0truein
 \hoffset=0.2truein
  \TagsOnRight

\nologo
 \NoRunningHeads
 \NoBlackBoxes

\def\ms{\medskip}

\def\ni{\noindent}

\def\p{\partial}

\def \D {\triangle}
\def\a{\alpha}
\def \e {\varepsilon}

\def \l {\lambda}

\def \o {\theta}

\def\W{\Omega}

\def\R{\Bbb{R} }
\def \li{\langle}
\def \ri{\rangle}
\def\Hess{\roman{Hess\/}}
\def\Exp{\roman{Exp\/}}
\def\gexp{\roman{gexp\/}}

\def\II{\roman{I\!I\/}}
\def\Int{\roman{Int\/}}
\def\dim{\roman{dim}}
\def\ang{\measuredangle}

\topmatter
\title An extension of Perelman's soul theorem for singular spaces
\endtitle

\author Jianguo Cao${}^1$, Bo Dai${}^2$ and Jiaqiang Mei${}^3$ \endauthor

\affil {University of Notre Dame, Peking University and Nanjing
University} \endaffil

\address Department of Mathematics, University of Notre Dame, Notre Dame, IN 46556, USA;
Department of Mathematics, Nanjing University, Nanjing 210093,
China.
\endaddress
\email cao.7\@nd.edu \endemail
\thanks ${}^1$The first author is supported in part by an NSF grant and Changjiang Scholarship of
China at Nanjing University. \endthanks

\address LMAM, School of Mathematical Sciences, Peking University, Beijing 100871, People's Republic of China
\endaddress
\email daibo\@math.pku.edu.cn \endemail
\thanks ${}^2$The second author is supported in part by NSFC Grants 10621061, 10701003. \endthanks

\address Department of Mathematics and Institute of Mathematical Science, Nanjing University,
Nanjing 210093, People's Republic of China
\endaddress
\email meijq\@nju.edu.cn \endemail
\thanks ${}^3$The third author is supported by Natural Science
Foundation of China (grant 10401015). \endthanks

\dedicatory In memory of Professor Xiao-Song Lin
\enddedicatory

\keywords Singular spaces, Perelman's soul theorem,
Cheeger-Gromoll convex exhaustion, generalized soul theory, angular excess functions
\endkeywords



\abstract In this paper, we study open complete metric spaces with
non-negative curvature. Among other things, we establish an
extension of Perelman's soul theorem for possibly singular spaces:
{\it ``Let $X$ be a complete, non-compact, finite dimensional
Alexandrov space with non-negative curvature. Suppose that $X$ has
no boundary and has positive curvature on a non-empty open subset.
Then  $X$ must be a contractible space"}. The proof of this result
uses the detailed analysis of concavity of distance functions and
Busemann functions on singular spaces with non-negative curvature.
We will introduce a family of angular excess functions to measure
convexity and extrinsic curvature of convex hypersurfaces in
singular spaces. We also derive a new comparison for trapezoids in
non-negatively curved spaces, which led to desired convexity
estimates for the proof of our new soul theorem.
\endabstract

\endtopmatter

\document
\baselineskip=.25in

\head \S 0. Introduction
\endhead

In this paper, we will study open complete and finite dimensional metric spaces with
non-negative curvature. A metric space $(X,d)$ is called a length space if any pair of
points $\{p,q\}$ in $X$ can be joined by a path of length equal to $d(p,q)$. A length-minimizing
path of unit speed is called a geodesic. A metric space $(X,d)$ is said to have non-negative
curvature if any geodesic triangle $\D_{\{l_1,l_2,l_3\}}$ of side lengths $\{l_1,l_2,l_3\}$ is
 ``fatter" than the comparison triangle $\D^*_{\{l_1,l_2,l_3\}}$ of the same lengths on the Euclidean plane. Similarly,
we can define the notion of curvature $\geq k$ for metric spaces with any real number $k\in\R$, see
\cite{BGP92}, \cite{BBI02}.

It is known that if a metric space $(X,d)$ has curvature bounded from below, then the topological
dimension of $X$ is equal to its Hausdorff dimension. Moreover, $\dim(X)$ must be an integer or
infinity (cf. \cite{BGP92}). Our main result of this paper is the following:

\proclaim{Theorem 0.1}  Let $X^n$ be a  complete and non-compact, $n$-dimensional metric space
of non-negative curvature. Suppose that $X^n$ has no boundary and has positive curvature  on a
metric ball $B_{\e_0}(x_0)$. Then $X^n$ must be contractible.
\endproclaim

When $X^n$ is a smooth open Riemannian manifold of non-negative curvature, our Theorem $0.1$ above
is related to the so-called Cheeger-Gromoll soul conjecture (cf. \cite{CG72}), which was successfully
solved by Perelman (cf. \cite{Per94a}).

However, there are examples of open metric spaces with positive curvature which are not topological
manifolds. For instance, let
$$M^n=\{(x_1,\dots, x_n,x_{n+1})\in\R^{n+1}\,|\,x_{n+1} =x_1^2 + \cdots + x_n^2 \}$$
and $X^n =M^n/ \Bbb{Z}_2$, where $\Bbb{Z}_2$ is a group generated by the involution
$$\align
\varphi: \R^{n+1} & \to \R^{n+1} \\
(x_1,\dots, x_n,x_{n+1}) & \mapsto  (-x_1,\dots, -x_n,x_{n+1}).
\endalign
$$
It is known that $X^n$ is a space of positive curvature. The space $\Sigma_0^{n-1}(X^n)$ of
unit tangent directions of $X^n$ at the origin $0$ is homeomorphic to the real projective
space $\R P^{n-1}$. Thus $X^n$ is not a manifold near the origin.

For smooth open Riemannian manifolds with non-negative curvature, Perelman (cf. \cite{Per94})
established a flat strip theorem and hence provided an affirmative solution to the
Cheeger-Gromoll soul conjecture. The proof of Perelman's flat strip theorem uses the
fact that, for a smooth Riemannian manifold $M^n$, its tangent space $T_x M^n$ at $x$ is
always isometric to $\R^n$. For an Alexandrov space $X^n$ with non-negative curvature, its
tangent cone $T_x X^n$ at a point $x$ is not necessarily isometric to $\R^n$. Therefore, a
different method is needed for the verification of our main theorem above. In next section,
we sketch our new approach by the study of convexity of Busemann functions on non-negatively
curved singular spaces.

\head \S 1. Outline of the proof of main theorem
\endhead

Our proof of Theorem $0.1$ is inspired by H.Wu's proof (cf. \cite{Wu79}, \cite{Wu87}) of
Gromoll-Meyer theorem (cf. \cite{GM69}).

\proclaim{Proposition 1.1}\rom(\cite{GM69}, \cite{Wu79}, \cite{Wu87}\rom) Suppose that $M^n$ is
a complete and non-compact smooth Riemannian manifold with positive sectional curvature, and
suppose that
$$h(x)=\lim_{t\to +\infty} [d(x, \p B_t (x_0)) - t ]$$
is a Busemann function, where $B_t (x_0) =\{ y\in M^n |d(x_0, y) <t \}$ is a metric ball of
radius $t$ centered at $x_0$. Then the function $(1-e^{-h})$ is a proper and strictly concave
function with a unique maximum point $\hat{p}\in M^n$. Consequently, $M^n$ is contractible to
a point $\hat{p}$ and hence $M^n$ is diffeomorphic to $\R^n$.
\endproclaim

Perelman (cf. \cite{Per91}) was able to derive a similar result for singular spaces as well. To
proceed, we need to recall the notion of $\lambda$-concavity introduced by Perelman (cf. \cite{Per94a})
for functions defined on singular spaces.

\proclaim{Definition 1.2}\rom($\lambda$-concave functions \cite{Per94}\rom) Let $f: X^n \to \R$
be a continuous function defined on an Alexandrov space with curvature $\geq -1$. We say that $f$
is $\lambda$-concave at $p$ (or $\Hess(f)\big|_p\leq\lambda$) in barrier sense if for all
quasi-geodesics $\sigma: (-\varepsilon, \varepsilon)\to X^n$ of unit speed with $\sigma(0)=p$, the
inequality
$$\frac{d^2}{dt^2}[f(\sigma(t))]\big|_{t=0}\leq\lambda$$
holds in barrier sense.
\endproclaim

\proclaim{Proposition 1.3}\rom(Perelman \cite{Per91} \rom Chapter 6) Let $X^n$ be a complete,
non-compact and $n$-dimensional metric space of positive curvature. Suppose that $X^n$ has no
boundary and
$$h(x)=\lim_{t\to +\infty} [d(x, \p B_t (x_0)) - t ]$$
is a Busemann function. Then the function $f(x)=\big(1-e^{-h(x)}\big)$ is a strictly concave
and a proper function with a unique maximum point $\hat{p}\in X^n$. Consequently, $X^n$ is
contractible to the maximum point $\hat{p}$ via the Sharafutdinov semi-flow
$$\frac{d^{+}\varphi}{dt} =\frac{\nabla h}{|\nabla h|^2}.$$
\endproclaim

We will review the Perelman-Sharafutdinov semi-gradient flow in upcoming sections.

We would like to say a few words about why we used $f=1-e^{-h}$ instead of the Busemann function
$h$. It is known (cf. \cite{Wu79}) that, for any $x^{*}\in X^n$, there is a geodesic ray
$\sigma: [0, +\infty)\to X^n$ of unit speed such that $\sigma(0)=x^{*}$ and
$$h(\sigma(t))=h(x^{*}) -t.$$
It follows that $h\circ\sigma$ is a linear function and hence $h$ can not be strictly concave
along the geodesic ray $\sigma$. Hence, it is reasonable to consider $f=1-e^{-h}$.

In our case, the singular space $X^n$ has positive curvature on a small ball $B_{\e_0}(x_0)$.
The function $f(x)=1-e^{-h(x)}$ is only weakly concave on the whole space $X^n$. The proof
of Perelman's result in Proposition $1.3$ also implies the following.

\proclaim{Proposition 1.4}\rom(Perelman \cite{Per94a} \rom) Let $X^n$ be a complete, non-compact
$n$-dimensional metric space of non-negative curvature. Suppose that $X^n$ has no boundary
and has positive curvature on a metric ball $B_{\e_0}(x_0)$ and
$$h(x)=\lim_{t\to +\infty} [d(x, \p B_t (x_0)) - t ]$$
as above. Then the function $f(x)=1-e^{-h(x)}$ is strictly concave on a small ball $B_{\e_0/4}(x_0)$.
\endproclaim

Inspired by Proposition $1.4$, we will take a close look on the concavity of $f(x)=1-e^{-h(x)}$
outside the small ball $B_{\e_0/4}(x_0)$. Following Cheeger-Gromoll (cf. \cite{CG72}), we consider
the convex sup-level sets:
$$\omega_c = h^{-1}([c, +\infty))$$
for all $c\in\R$. it is known (cf. \cite{Per94a} or \cite{CMD09}) that $\Omega_c$ is a totally
convex subset of $X^n$. Moreover, its boundary $\p\Omega_c$ has strictly convex portion
$(\p\Omega_c)\cap B_{\e_0/4}(x_0)$ when $c\in[\hat{c}_0-\frac{\e_0}{4}, ~\hat{c}_0+\frac{\e_0}{4}]$
and $\hat{c}_0=h(x_0)$. It is also known that each $\Omega_c$ is compact. Thus, we may assume
$$c_0 = \max_{x\in X^n} \{h(x)\} < +\infty. \tag 1.1$$
Hence, $\Omega_c=h^{-1}([c, +\infty))=h^{-1}([c, c_0])$ for $c\leq c_0$.

If $c_0=\hat{c}_0$ and $f(x)=1-e^{-h(x)}$ is strictly concave at $x_0$ then $x_0$ is the unique
maximum of $f$ and $h$; hence $X^n$ is contractible. Thus, we may assume that $\hat{c}_0=h(x_0)<c_0$.

There are three possibilities for the maximum set $\Omega_{c_0}=A_0$ of $h$:

{\bf Case 1}. $\Omega_{c_0}$ is a convex and compact subset without boundary and $\dim(\Omega_{c_0})\geq 1$.
In this case, $\Omega_{c_0}$ remains an Alexandrov space of non-negative curvature.

{\bf Case 2}. $\dim(\Omega_{c_0})=0$ and $X^n$ is contractible.

{\bf Case 3}. $\dim(\Omega_{c_0})\geq 1$ but $\Omega_{c_0}$ is a convex subset with non-empty
boundary. In this case, we let
$$\Omega_{c_0+s} =\{x\in\Omega_{c_0}\,|\,d(x, \p\Omega_{c_0})\} \geq x\}$$
and consider the distance function
$$r_{\p\Omega_{c_0}}(x)=d(x, \p\Omega_{c_0})  \tag 1.2$$
for $x\in\Omega_{c_0}$. Since $\Omega_{c_0}=A_0$ is compact, the distance function $r_{\p\Omega_{c_0}}$
from the boundary has a maximum value on $A_0$, say
$$l_1=\max\{r_{\p\Omega_{c_0}}(x)\,|\,x\in\Omega_{c_0}\}. $$
There are also three possibilities for the maximum subset $A_1=\Omega_{c_0+l_1}=r^{-1}_{\p\Omega_{c_0}}(l_1)$
as well.

Since we have
$$\dim(X^n)>\dim(A_0)>\dim(A_1)>\cdots,$$
repeating above operations at most $n$ times, we end up either {\bf Case 1} or {\bf Case 2}.
Under the assumption that $X^n$ has positive curvature on $B_{\e_0}(x_0)$, we will make efforts
to rule out {\bf Case 1} above.

\proclaim{Theorem 1.5} Let $X^n$ be a complete, non-compact $n$-dimensional metric space of
non-negative curvature and
$$h(x)=\lim_{t\to +\infty} [d(x, \p B_t (x_0)) - t ]$$
as above. Suppose that $X^n$ has no boundary and has positive curvature on a metric ball
$B_{\e_0}(x_0)$. Then the maximum subset $A_0=h^{-1}(c_0)$ of $h$ must be either a point
or $A_0$ has a (strictly convex) boundary point $y_0\in\p A_0$ in the sense defined below.
\endproclaim

The proof of Theorem $1.5$ will be given in upcoming sections.

For a convex subset $A\subset X^n$, there is sufficient condition for the subset $A$ to have
a boundary point.

\proclaim{Definition 1.6} Let $X^n$ be an Alexandrov space of
curvature $\geq -1$.

$(1.6.l)$ If, for any quasi-geodesic $\sigma: [a, b]\to X^n$with ending points
$\{\sigma(a), \sigma(b)\}\subset A$, quasi-geodesic segment $\sigma([a, b])\subset A$, then
$A$ is called a totally convex subset of $X^n$.

$(1.6.2)$ Suppose that $A$ is a totally convex subset $A$ of $X^n$, $p\in A$ and $\vec{u}\in T_p(X^n)$
is a unit tangent direction of $X^n$ at $p$; suppose that
$$\ang_p(\vec{u}, x) \geq \frac{\pi}{2} \tag 1.3$$
for any $x\in A -\{p\}$. Then $\vec{u}$ is called at least normal to $A$ at $p$.

$(1.6.3)$ If $\vec{u}\in T_p(X^n)$ and $\|\vec{u}\|=1$ then the sub-space
$$H^{+}_{\vec{u}} =\{y\in X^n-\{p\}\,|\,\ang_p(\vec{u}, y)\}>\frac{\pi}{2}$$
is called an open half sub-space relative to $\vec{u}$.

$(1.6.4)$ Suppose that $\vec{u}$ is an at least normal vector to a totally convex subset $\Omega$
at $p$, and suppose that
$$\theta_{p,\vec{u}}^{\Omega}(r)=\inf\{\ang_p(\vec{u},x)|x\in\Omega, d(x,p)\geq r\} -\frac{\pi}{2}>0,  \tag 1.4$$
for all sufficiently small $r>0$. Then the point $p$ is called a strictly convex boundary point
of $\Omega$.
\endproclaim

When $X^n$ has positive curvature on $B_{\e_0}(x_0)$, we already point out that the function
$f(x)=1-e^{-h(x)}$ is strictly concave on a smaller ball $B_{\e_0/4}(x_0)$. Moreover, we have
the following refined estimate.

\proclaim{Proposition 1.7} Let $X^n$, $h(x)$ and $\Omega_c=h^{-1}([-c, c_0])$ be as above. If
$X^n$ has positive curvature $\geq k_0>0$ on $B_{\e_0}(x_0)$, then for each $p\in B_{\e_0/4}(x_0)$
there exists $a>0$ such that
$$\theta_{p,\vec{u}}^{\Omega_c}(r) \geq ar>0  \tag 1.5$$
for some $\vec{u}\in T_p(X^n)$ and sufficiently small $r>0$, where $c=h(p)$.
\endproclaim

Let us now return to the maximum subset $A_0=\Omega_{c_0}=h^{-1}(c_0)$. For interior points
of $A_0$, we have the following observation.

\proclaim{Proposition 1.8} \rom(Compare \cite{CG09} \S 2.2\rom) Let $X^n$, $h$, $\{\Omega_c\}$
and $\theta_{p,\vec{u}}^{\Omega}(r)$ be as above. Suppose that $\dim(A_0)\geq 1$ and $p$ is an
interior point of $A_0$. Then
$$\ang_p(\vec{u}, y) = \frac{\pi}{2}  \tag 1.6$$
for any at least normal direction $\vec{u}$ to $A_0$ at $p$ and any $y\in A-\{p\}$. Consequently,
$$\theta_{p,\vec{u}}^{A_0}\equiv 0  \tag 1.7$$
for all $r>0$.
\endproclaim

\vskip.3in

\epsfxsize=5.5cm             
\centerline{\epsfbox{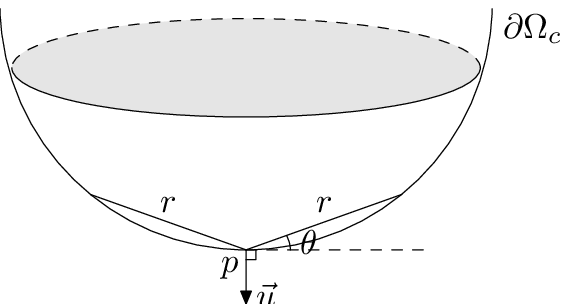}}   
\botcaption{Figure 1. Strictly convexity and angular excess}
\endcaption

In order to establish Theorem $1.5$ and hence Theorem $0.1$, by $(1.7)$ it is more desirable
to show that the inequality $(1.5)$
$$\theta_{p,\vec{u}}^{\Omega_c}(r) \geq ar>0 $$
for some points $p\in A_0$ and $\Omega_c=A_0$.

Thus, we are led to study the inequality $(1.5)$ along a
Perelman-Sharafutdinov curve
$$\frac{d^{+}\varphi}{dt} = \frac{\nabla h}{|\nabla h|^2}.$$

It was shown by H.Wu (cf. \cite{Wu79}) that our Busemann function $h$ is indeed a distance
from an appropriate subset. In fact, one can show that
$$h(x)=c+d(c, \p\Omega_c)=c+r_{\p\Omega_c}(x),  \tag 1.8$$
for all $x\in\Omega_c=h^{-1}([c, +\infty))$. Therefore, the semi-gradient semi flow of $h$
is actually semi-gradient flow for distance functions. The evolution induced by Perelman-Sharafutdinov
semi-flow
$$\frac{d^{+}\varphi}{dt} = \frac{\nabla h}{|\nabla h|^2} = \frac{\nabla r_{\p\Omega_c}}{|\nabla r_{\p\Omega_c}|^2}  \tag1.9$$
is indeed an one-sided equidistance evolution. To our surprise, the angular estimate of
$$\theta_{p,\vec{u}}^{\Omega_c}(r) \geq ar>0$$
is preserved by Perelman-Sharafutdinov equi-distance evolution $\{\p\Omega_c\}$ in a space $X^n$
of non-negative curvature.

\proclaim{Theorem 1.9} \rom(Angular estimates under parallel translation \rom) Let $X^n$, $h(x)$,
$\Omega_c=h^{-1}([c, +\infty))$ and $\theta_{p,\vec{u}}^{\Omega_c}(r)$ be as above. Suppose
that $X^n$ has non-negative curvature, $\varphi: [a, b]\to X^n$ is a Perelman-Sharafutdinov
curve
$$\frac{d^{+}\varphi}{dt} = \frac{\nabla h}{|\nabla h|^2}$$
with $c_i=h(\varphi(t_i))$, $t_1\leq t_2$, $\vec{u}$ is the left-derivative of $\varphi$ at
$t$, and suppose that $\theta_t(r)=\theta_{\varphi(t),\vec{u}(t)}^{\Omega_{h(\varphi(t))}}(r)$.
Then
$$\theta_{t_2}(r) \geq \theta_{t_1}(r)  \tag 1.10$$
for any $t_1\leq t_2$. Equivalently,
$$\frac{\p\theta_t(r)}{\p t} \geq  0  \tag 1.11$$
holds.
\endproclaim

\vskip.3in

\epsfxsize=8cm             
\centerline{\epsfbox{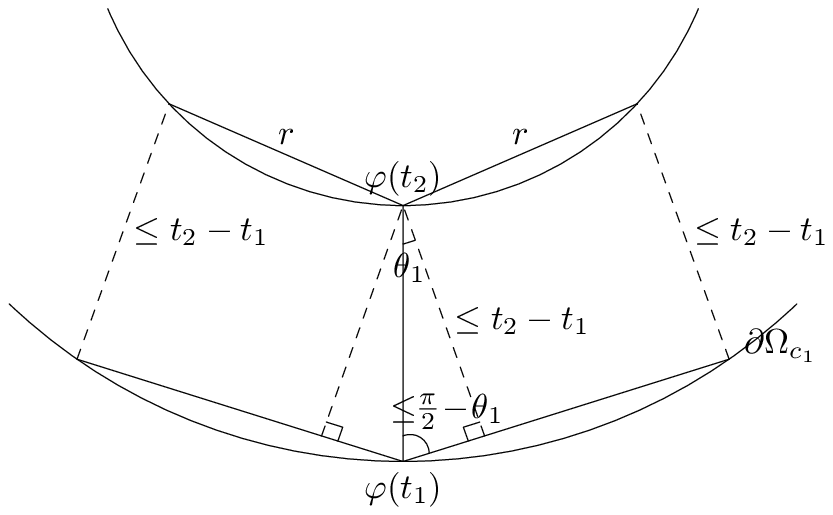}}   
\botcaption{Figure 2. Angular estimates under equidistance evolution}
\endcaption

If $X^n$ is a smooth Riemannian manifold of non-negative curvature
and if $\p\Omega_{c_i}$ is a smooth convex hypersurface then the
inequality
$$\frac{\p\theta_t(r)}{\p t} \geq  0  $$
is related to the classical Riccati equation of the smooth second fundamental form of
$\II_{\p\Omega_{h(\varphi(t))}}=\II_t$:
$$\II'_t + \II^2_t + R =0,  \tag 1.12$$
where $R(V, Y)Z = -\nabla_V\nabla_Y Z = \nabla_Y\nabla_V Z + \nabla_{[V, Y]} Z$ is the
curvature tensor of the smooth Riemannian manifold $X^n$.

For instance, let $\Omega\subset\R^2$ be a domain given by
$$\Omega=\{(x_1, x_2)\,|\,x_2<-x_1^2\}.$$
If we choose $p=(0,0)$ and $\vec{u}=(0, 1)$, then $\theta_{p, \vec{u}}^{\Omega}(r)\sim\tan^{-1}r\sim r$
for all sufficiently small $r$.

\vskip.3in

\epsfxsize=6cm             
\centerline{\epsfbox{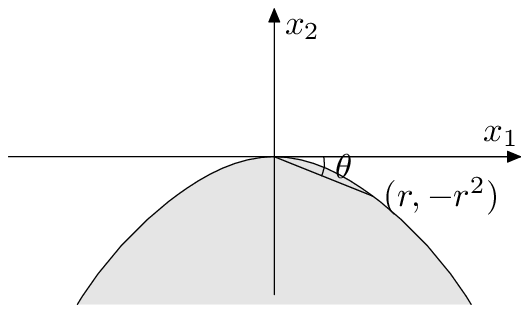}}   
\botcaption{Figure 3.}
\endcaption

\proclaim{Proposition 1.10} Let $\Omega\subset M^n$ be a totally
convex subdomain with smooth boundary $\p\Omega$ in a smooth
Riemannian manifold $M^n$. Then the second fundamental form
$\II_{\p\Omega}$ of $\p\Omega$ with respect to the outward unit
normal vector $\vec{u}$ at $p$ satisfies
$$\II_{\p\Omega}^{\vec{u}}(\vec{w}, \vec{w}) = -\li\nabla_{\vec{w}}\vec{w}, \vec{u}\ri
\geq 2\lambda_1\|\vec{w}\|^2   \tag 1.13$$
for some $\lambda_1>0$ if and only if
$$\lim_{r\to 0} \frac{\theta_{p,\vec{u}}^{\Omega}(r)}{r} \geq \lambda_2>0  \tag 1.14$$
for some $\lambda_2>0$.
\endproclaim

In a recent paper \cite{AB09}, Alexander and Bishop used the length-excess function for chords
in $\Omega$ to define extrinsic curvature of $\p\Omega$ in Alexandrov spaces. We will discuss
the relation between our angular excess functions $\liminf\limits_{r\to 0}\frac{\theta_{p,\vec{u}}^{\Omega}(r)}{r}$
and Alexander-Bishop's extrinsic curvature in upcoming sections.

Among other things, we will use the following trapezoid comparison
theorem to verify Theorem $1.9$ above.

\vskip.3in

\epsfxsize=6.5cm             
\centerline{\epsfbox{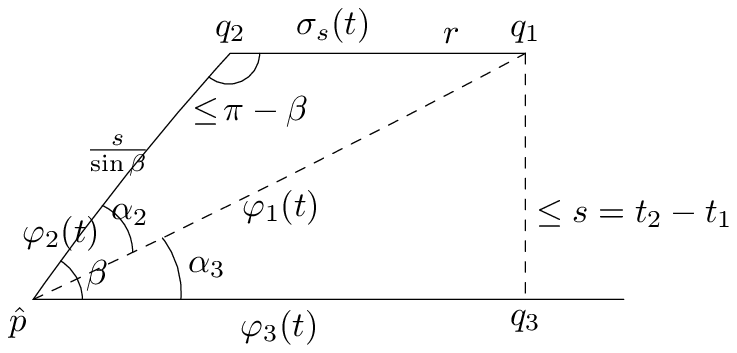}}   
\botcaption{Figure 4. Trapezoid comparison}
\endcaption

\proclaim{Theorem 1.11} \rom(Trapezoid comparison theorem \rom) Let $X^n$ be a complete
Alexandrov space of non-negative curvature as above, and let $\{\varphi_1, \varphi_2, \varphi_3\}$
be three geodesic segments with the same initial point $\hat{p}$. Suppose that the three
initial directions $\{\varphi_1'(0), \varphi_2'(0), \varphi_3'(0)\}$ are co-planar in the sense
$$\align
& \ang_{\hat{p}}(\varphi_2'(0), \varphi_1'(0)) + \ang_{\hat{p}}(\varphi_1'(0), \varphi_3'(0)) \\
& = \alpha_2 + \alpha_3 =\beta \\
& = \ang_{\hat{p}}(\varphi_2'(0), \varphi_3'(0)),   \tag 1.15
\endalign
$$
$l=d(\hat{p}, q_1)$, $\varphi_3$ is a possible quasi-geodesic, $d(\hat{p}, q_2)=\frac{s}{\sin\beta}$
and $\sigma_3: [0, r]\to X^n$ is a geodesic segment from $q_2$ to $q_1$ with
$$\ang_{q_2}\big(\sigma_s'(0), \frac{d^{-}\varphi}{dt}(\frac{s}{\sin\beta})\big)\leq \pi-\beta  \tag1.16$$
as in Figure 4. Then the inequality
$$d(\sigma_s(r), \varphi_3(\R)) \leq d(\sigma_s(r), \varphi_3(l\cos\alpha_3))\leq s  \tag 1.17$$
holds.
\endproclaim

We provide the detailed proofs of results stated above in upcoming sections.

\head \S 2. Non-negative curvature and weak concavity of Busemann functions
\endhead

In this section, we will prove Proposition $1.8$ which is related to the weak concavity of
Busemann functions defined on open complete spaces of non-negative curvature.

We will use the same notation as in Section $1$. Let $h: X^n \to \R$ be a Busemann function
given by
$$h(x)=\lim_{t\to +\infty} [d(x, \p B_t (x_0)) - t ]$$
and let $\Omega_c=h^{-1}([c, +\infty))$. In order to derive the
desired angular estimates, we recall the totally convex subsets so
that we can construct conic-like barrier hypersurfaces.

\proclaim{Proposition 2.1} Let $X^n$ be an open, complete and
$n$-dimensional Alexandrov space with non-negative curvature and
let $h(x)$ and $\Omega_c = h^{-1}([c, +\infty))$ be as above. Suppose that $\sigma: [a,b] \to
X^n$ is a quasi-geodesic of unit speed with $\{\sigma(a), \sigma(b)\} \subset \Omega_c$. Then
$$\sigma([a,b]) \subset \Omega_c \tag 2.1 $$
and $\Omega_c$ is a totally convex subset of $X^n$.
\endproclaim

\demo{Proof} The proof of Proposition $2.1$ for smooth Riemannian manifold was proved by Cheeger-Gromoll
(cf. \cite{CG72}). For singular spaces, Perelman (cf. \cite{Per91} Chapter 6) also established
the similar result (see also Section 2 of \cite{CaG10}). For convenience of readers, we present
another simple proof here. By an equivalent definition of non-negative curvature, for any distance
function $r(x) = d(x, \hat y)$ from a given point $\hat y$, we have
$$
            \Hess(  \frac 12 r^2) \le I \tag2.2
$$
which means that, for any quasi-geodesic $\sigma $ of unit speed, we have
$$
  \frac{d^2}{ds^2} \big( \frac 12 [r(\sigma(s))]^2 \big)  \le 1. \tag2.3
$$
It follows that
$$
    \Hess(r) \le \frac{1}{r} I. \tag2.4
$$
Recall that $h(x) = \lim_{t \to \infty} [ d(x, \partial B_t(\hat x)) - t]$. By (2.4) we have
$$
        \Hess(h(x)) = \lim_{ t \to \infty} [  \Hess(  \frac 12 [d(x, \partial B_t(\hat x))]^2)] \le
         \lim_{ t \to \infty}[ \frac{1}{d(x, \partial B_t(\hat x))}] =0.
$$
Hence $h$ is a concave function. It follows that
$$
   h(\sigma(t)) \ge \min\{h (\sigma(a)), h (\sigma(b)) \}
$$
for all $t \in [a, b]$. Therefore, the sup-level set $\Omega_c = h^{-1}([c, +\infty))$ is a
totally convex subset of $X^n$.
\qed
\enddemo

We also need to show that the Busemann function $h$ is proper.

\proclaim{Proposition 2.2} \rom(Compare \cite{Wu79} \rom) Let $X^n$ be an open, connected
and complete $n$-dimensional space of non-negative curvature and
$h(x) = \lim\limits_{t \to \infty} [ d(x, \partial B_t(\hat x)) - t]$ as above. The for each $c\in\R$,
the sup-level set $\Omega_c=h^{-1}([c, +\infty))$ is compact.
\endproclaim

\demo{Proof} Let $\hat{c}=\min\{h(\hat{x}), c\}=\min\{0, c\}$. It
is sufficient to verify that $\Omega_{\hat{c}}$ is compact.
Suppose to the contrary, $\Omega_{\hat{c}}$ were non-compact.
There would be un-bounded sequence
$\{q_i\}\subset\Omega_{\hat{c}}$ such that
$$l_j=d(\hat{x}, q_j) \to +\infty  \tag2.5$$
as $j\to+\infty$. Let $\sigma_j: [0, l_j]\to X^n$ be a
length-minimizing geodesic from $\hat{x}$ to $q_j$. By passing to
a subsequence $\{y_{j_k}\}$ of $\{y_j\}$ if necessary, we may
assume that $\sigma'_j(0)\to\vec{v}_{\infty}$ and
$\sigma_j\to\sigma_{\infty}$ as $j\to+\infty$, where
$\sigma_{\infty}: [0, +\infty) \to X^n$ is a geodesic ray from
$\hat{x}$. H.Wu (cf. \cite{Wu79}) observed that
$$h(\sigma_{\infty}(t))=h(\hat{x})-t=-t\to-\infty<\hat{c}  \tag2.6$$
as $t\to+\infty$. However, we already proved that $\Omega_{\hat{c}}$ is totally convex, and
hence
$$\sigma_{\infty}([0, +\infty)) \subset \Omega_{\hat{c}},   \tag 2.7$$
which implies
$$h(\sigma_{\infty}(t)) \geq \hat{c},    \tag2.8$$
a contradiction to $(2.6)$.
\qed\enddemo

By Proposition $2.2$, we see that our Busemann function $h$ has a finite maximum value $a_0$
where
$$a_0=\max_{x\in X^n} \{h(x)\}.  \tag2.9$$
We consider the maximum set $A_0=h^{-1}(a_0)$ of $h$. When $\dim(A_0)\geq 1$, we would like to
address the ``weak concavity property" of $A_0=h^{-1}(a_0)$ at its interior points.

\proclaim{Proposition 2.3} \rom(Compare Theorem 2.11(2) in \cite{CaG10} \rom) Let $X^n$ be an
open complete $n$-dimensional space of non-negative curvature,
$h(x)=\lim\limits_{t\to +\infty} [d(x, \p B_t(x_0))-t]$ and $a_0=\max \{ h(x)|x\in X^n\}$.
Suppose that $A_0=h^{-1}(a_0)$ has positive dimension $\dim(A_0)\geq 1$, $x\in \roman{int}(A)$
is an interior point of $A_0$, and $\vec v$ is a unit vector at least normal to $A_0$ at $x$.
Then
$$\ang_x (\vec v,\vec w) =\frac\pi{2} \tag 2.10$$
holds for all $\vec w \in T_x(A)$.
\endproclaim

\demo{Proof} {\bf Case 1.} $\dim(A_0)=1$. Since $A_0$ is a totally convex subset of $X^n$, $A_0$
is either a closed geodesic or a length-minimizing geodesic segment. By our assumption, $x$ is
an interior point of $A_0$, there is a geodesic $\sigma: (-\e,\e) \to A_0 \subset X^n$ such that
$\sigma(0)=x$. It is known (cf. \cite{BGP92}) that the tangent space $T_x(X^n)$ has an isometric
splitting
$$ T_x(X^n) =Y^{n-1} \times T_x(A) =Y^{n-1} \times \R. \tag 2.11$$
It follows that, for any $\vec v\in \Sigma_x(X^n)$, we have
$$ \min\{ \ang_x (\vec v, \sigma'(0)),  \ang_x (\vec v, -\sigma'(0)) \} \leq \frac\pi{2}. \tag 2.12 $$
If $\vec v$ is at least normal to $A_0$ at $x$, then by $(2.12)$, we must have
$$ \ang_{x} (\vec v, \pm \sigma'(0)) =\frac\pi{2}. $$
and $\vec v \in Y^{n-1} \perp T_x(A_0)$.

{\bf Case 2.}  When $\dim(A_0) >1$, our proof becomes more
involved. If $Y$ is an Alexandrov space with curvature $\geq c$
and $y\in Y$, then we let $\Sigma_y(Y)$ denote the space of unit
tangent directions of $Y$ at $y$. It is well-known that
$\Sigma_y(Y)$ has curvature $\geq 1$. Let $A_0'=\Sigma_x(A_0)$,
where $A_0=\{y\in X \ |\,h(y)=\max_{z\in X}\{h(z)\}\}$ is the
maximum set of Busemann function $h$ defined on the open space
with non-negative curvature. Using triangle comparison theorem for
spaces $\Sigma$ with curvature $\geq 1$ and a result of
Perelman-Petrunin on quasi-geodesics, we will show that if $x$ is
an interior point of $A_0$, then
$$
\max\{\ang_x(\vec{w}, A_0')\,\big|\,\vec{w}\in\Sigma_x(X), ~\ang_x(\vec{w}, A_0')=d_{\Sigma}(x, A_0')\}\leq\frac{\pi}{2}.  \tag 2.13
$$
Our strategy to establish $(2.13)$ goes as follows (compare with the proof of Theorem 2.11(2) in
\cite{CaG10}).

\smallskip
\noindent
{\bf Claim A.} {\it Let $A_0$ be a totally convex and the maximum set of a Busemann function $h$
as above. Suppose that $x$ is an interior point of $A_0$, $A_0'=\Sigma_x(A_0)$ and
$\dim(A_0) = k \ge 2$. Then the following is true:

\smallskip
\noindent
(A.1) For each $\vec{u}\in\Sigma_{\xi}(A_0')$, there is a quasi-geodesic $\sigma_v: (-\e, \e)\to A_0'$
such that $\sigma_v(0)=\vec{\xi}$ and $\sigma'_v(0)=\vec{u}$.

\smallskip
\noindent
(A.2) $d_{\Sigma}(\vec{w}, A_0')\leq \frac{\pi}{2}$ for all $\vec{w}\in\Sigma_x(X)$.}

\smallskip

We first verify Assertion (A.1). Recall that $A_0$ is totally convex in $X$. Therefore,
$A_0'=\Sigma_x(A_0)$ is totally convex in $\Sigma_x(X)$. Since $x$ is an interior point of $A_0$,
the subspace $A_0'=\Sigma_x(A_0)$ has no boundary in $\Sigma_x(X)$ with $\dim(A_0')=k-1\geq 1$.
Hence, $A_0'$ is a compact, totally convex subspace without boundary in $\Sigma_x(X)$. For each
$\xi\in A_0'$, we let
$$ \align
\rho_{\xi} : \Sigma_x(X) & \longrightarrow \R \\
 \vec{w} & \longrightarrow \ang_x(\xi, w)=d_{\Sigma}(\xi, w).
\endalign
$$
Since $A_0'$ is a totally convex subspace without boundary in $\Sigma_x(X)$, the semi-gradient
curves $\varphi: [0, l)\rightarrow \Sigma_x(X)$ of $\frac{d^{+}\varphi}{dt}=\nabla\rho_{\xi}$
with $\varphi(0)=\xi$ and $\varphi'(0)=\vec{u}$ will remain in $A_0'$. Moreover, the gradient
exponential map at $\xi$, $\gexp_\xi: T_\xi(\Sigma_x(X)) \to \Sigma_x(X)$,  has the property
$\gexp_\xi(T_\xi(A_0')) \subset A_0'$, because $A_0'$ is totally convex in $\Sigma_x(X)$ and
because $(T_\xi(\Sigma_x(X)), o_\xi) = \lim\limits_{\lambda \to \infty}(\lambda \Sigma_x(X), \xi)$.
It follows from  a theorem of Perelman-Petrunin (cf. \cite{Petr07}) that quasi-geodesics can
be approximated by broken gradient exponential curves, (cf. Appendix of \cite{Petr07}). Thus,
by the construction of quasi-geodesics given by Perelman-Petrunin (cf. Appendix of \cite{Petr07}),
we see that any quasi-geodesic
$$\sigma: [0, l] \rightarrow \Sigma_x(X)$$
with $\sigma(0)=\xi$ and $\sigma'(0)=\vec{u}\in\Sigma_{\xi}(A_0')$ will stay in $A_0'\subset\Sigma_x(X)$.
This completes the proof of our Assertion (A.1).

For Assertion (A.2), we use a triangle comparison theorem for the space $\Sigma_x(X)$. Let
$\varphi: [0, l]\rightarrow \Sigma_x(X')$ be a length-minimizing geodesic segment from $A_0'$ to
$\vec{w}\in\Sigma_x(X)$ of unit speed. Suppose that $\varphi(0)=\xi\in A_0'$, $\varphi(l)=w$ and
$$d_{\Sigma}(A_0', w) =d(\xi, w)=l.$$
We will use comparison theorem to show that
$$d_{\Sigma}(A_0', w) = d(\xi, w)=l\leq\frac{\pi}{2}. \tag 2.14$$

Because $\dim(A_0')=k-1\geq 1$, we can choose a quasi-geodesic
$\sigma: (-\varepsilon, \varepsilon)\rightarrow A_0'\subset\Sigma_x(X)$ such that $\sigma(0)=\xi$.
Let $\sigma_{\pm}'(0)=\vec{u}_{\pm}$ be the left (or right) derivative of $\sigma$ at $s=0$.
Perelman and Petrunin (cf. \cite{Petr07}) showed that
$$\langle \vec{u}_{+}, \varphi'(0)\rangle + \langle \vec{u}_-, \varphi'(0)\rangle \geq 0.$$
It follows that
$$\alpha=\min\{\ang_{\xi}(\vec{u}_{+}, \varphi'(0)), ~\ang_{\xi}(\vec{u}_{-}, \varphi'(0))\} \leq \frac{\pi}{2}.$$
We may assume that $\vec{u}_{+}=\sigma'(0)$ has the property
$$\alpha=\ang_{\xi}(\vec{u}_{+}, \varphi'(0))=\ang_{\xi}(\sigma_{+}'(0), \varphi'(0)) \leq\frac{\pi}{2}.$$
If $l\leq\frac{\pi}{2}$, we are done. If $l>\frac{\pi}{2}$, then we will get a contradiction as follows.
For geodesic hinge $\{\sigma, \varphi\}$ with the vertex $\xi$ and angle $\alpha\leq\frac{\pi}{2}$,
since $\Sigma_x(X)$ has curvature $\geq 1$, we have
$$d_{\Sigma}(\varphi(\frac{\pi}{2}), \sigma(\varepsilon)) \leq \frac{\pi}{2}. \tag 2.15 $$
Recall that $\varphi: [0, l]\rightarrow \Sigma_x(X)$ is a length-minimizing geodesic from $A_0'$,
we also have
$$d_{\Sigma}(\varphi(\frac{\pi}{2}), \sigma(\varepsilon)) \geq d_{\Sigma}(\varphi(\frac{\pi}{2}), A_0')\geq \frac{\pi}{2}.$$
Combing with $(2.15)$, we see that
$$d_{\Sigma}(\varphi(\frac{\pi}{2}), \sigma(\varepsilon)) =\frac{\pi}{2}.$$
Therefore, we have two distinct points $\{\xi, \sigma(\varepsilon)\}\subset A_0'$ such that
$$d(\xi, \varphi(\frac{\pi}{2})) =d(\sigma(\varepsilon), \varphi(\frac{\pi}{2}))
  =d(A_0', \varphi(\frac{\pi}{2}))=\frac{\pi}{2}.$$
Thus, $\varphi|_{[0, \frac{\pi}{2}+\varepsilon']}$ is no longer a length-minimizing geodesic from
$A_0'$ for sufficiently small $\varepsilon'$ where $0 < \varepsilon' \le [l- \frac{\pi}{2}]$, a
contradiction. Thus we established
$$l=d(A_0', w)\leq \frac{\pi}{2}$$
for all $w\in\Sigma_x(X)$. This completes the proof of our Assertion (A.2) as well as Proposition $2.3$.
\qed
\enddemo

In next section, we discuss the relation between strictly concavity of $h$ and positive curvature
on the small ball $B_{\e_0}(\hat{x})$.

\head \S 3. Strictly positive curvature and strong concavity of Busemann function on a small region
\endhead

By our discussions in \S 2 above, we see that if our Busemann function $h$ is partially strong
concave on a portion of the maximum set $A_0=h^{-1}(a_0)$, then either $A_0$ is a single point
set or $A_0$ has non-empty boundary.

In order to establish the partial strong concavity of $h$ on a portion of the maximum set $A_0$,
we begin with a small ball $B_{\e_0/4}(x)$ where the curvature of $X^n$ is strictly positive.

\proclaim{Proposition 3.1} Let $X^n$ be an open complete Alexandrov space of non-negative curvature.
Suppose that $X^n$ has curvature $\geq 1$ on $B_{2\e_0}(p_0)$,
$$h(x) = \lim\limits_{t \to \infty} [ d(x, \partial B_t(\hat x)) - t]$$
and $c_0=h(\hat{x})$. Then, for any $p\in B_{\e_0}(\hat{x})$, there exists $\vec{u}$ with
$$\theta_{p,\vec{u}}^{\Omega_{h(p)}}(r)\geq \frac{\e_0}{8}r  \tag 3.1$$
for sufficiently small $0<r\leq r_0$.
\endproclaim

\demo{Proof} We will use a result of Petrunin (\cite{Petr07}) to derive the desired estimate.

By our discussion,
$$ h(x)=(c_0 -2\e_0) +d(x, \p\Omega_{c_0 -2\e_0}) $$
for $x\in B_{2\e_0} (p_0)$. We choose $\vec v\in\Uparrow_{p_0}^{\p \Omega_{c_0 -2\e_0}}$ and
$\sigma: [0,2\e_0] \to X^n$ is a length-minimizing geodesic from $p_0$ to $\p\Omega_{c_0 -2\e_0}$
with $\sigma'(0) =\vec v$. Let $q_0 =\sigma(\e_0)$, $\vec w \in \Sigma_{p_0} (X^n)$ such that
$\ang_{p_0} (\vec v, \vec w) =\frac\pi{2}$, and $\psi: [0,r_0] \to X^n$ be a quasi-geodesic
from $p_0$ with $\psi'(0)=\vec w$. Petrunin (\cite{Petr07}) proved that $d(\psi(r), \p\Omega_{c_0-\e})$
is bounded above by the $2$-dimensional model case as follows. Let
$$S^2=\{ (x_1,x_2,x_3)\in\R^3| x_1^2 +x_2^2 +x_3^2 =1\}, ~~S^1 =\{ (x_1,x_2,0)| x_1^2 +x_2^2 =1 \},$$
$p_0^* =(\cos\e_0,0,\sin\e_0)$ and $q_0^* =(1,0,0)$. The spherical geodesic
$$ \psi^*(r) =(\cos r)p_0 +(\sin r)z_0 =\big((\cos\e_0)\cos r, \sin r, (\sin\e_0)\cos r\big).$$
A calculation shows that
$$ \align
d(\psi^*(r), S^1) &=\tan^{-1} \Big( \frac{(\sin\e_0)\cos r}{\sqrt{(\cos \varepsilon_0)^2 (\cos r)^2 +\sin^2 r}} \Big) \\
&\leq \e_0\cos r \leq \e_0(1-\frac{r^2}{2} +\frac{r^4}{24}). \tag 3.2
\endalign
$$
Petrunin (\cite{Petr07}) showed that
$$ d(\psi(r), \p\Omega_{c_0-\e_0}) \leq d_{S^2} (\psi^*(r), S^1) +o(r^2)
   \leq \e_0  (1-\frac{r^2}{2} +\frac{r^4}{24}) +o(r^2),  \tag 3.3$$
where $\lim\limits_{r\to 0} \frac{o(r^2)}{r^2} =0$. Therefore, for sufficiently small $r< r_0$,
we have
$$ d(\psi(r), \p\Omega_{c_0-\e_0}) \leq \e_0  (1-\frac{r^2}{4}). \tag 3.4$$

It follows that
$$h(\psi(r)) \le c_0 - \frac{\varepsilon_0 r^2}{4}. $$
Thus for each $y\in [\Omega_{c_0} -B_r(p_0)]$ with $d(y, p_0) = r$, we have
$d(y,  \psi(r)) \ge \frac{\varepsilon_0 r^2}{4}$.  By comparison theorem for spaces
with curvature $\ge 0$, we see that
$$ \ang_{p_0} (y, \psi'(0))  \ge \frac{\varepsilon_0 r^2}{4r} \ge \frac{\e_0}{8} r. \tag 3.5$$
It follows that
$$ \o^{\Omega_{c_0}}_{p_0,\vec v} (r) \geq \frac{\e_0}{8} r \tag 3.6$$
for sufficiently small $r<r_0$.  \qed
\enddemo

In next section, we will discuss a sufficient condition so that the inequality
$$\theta_{p,\vec{u}}^{\Omega_{h(p)}}(r) > c_0 r$$
holds for some points $p\in A_0=h^{-1}(a_0)$.

\head \S 4. Preserving concavity of  hypersurfaces under equi-distance evolution
\endhead

In previous section, we showed that if $X^n$ has positive curvature on $B_{2\e_0}(\hat{x})$
then $B_{\e_0}(\hat{x})\cap\p\Omega_{h(x)}$ is strictly concave. Our goal of this section is
to show that the strictly concavity property of $\{\p\Omega_{h(\varphi(t))}\}$ is preserved
along equi-distance evolution. More precisely, if $\varphi: [a, b]\to X^n$ is a Perelman-Sharafutdinov
curve for a distance function (or a Busemann function) and if $\p\Omega_{h(\varphi(0))}$ is
strictly concave at $\varphi(0)$ then $\p\Omega_{h(\varphi(t))}$ remains strictly concave at
$\varphi(t)$ for all $t\geq 0$. Let us discuss several examples to motivate Theorem $1.9$.

{\bf Example 4.1} (a) We consider the following domain
$$\Omega_0=\{(x_1,x_2)\in\R^2\,|\,|x_1|\leq 4, ~|x_2|\leq 1, ~\text{or}
  ~|x_1|\geq 4, ~\frac{1}{4}(x_1\pm 4)^2+x_2^2\leq 1\},$$
see Figure 5.

\vskip.3in

\epsfxsize=11cm             
\centerline{\epsfbox{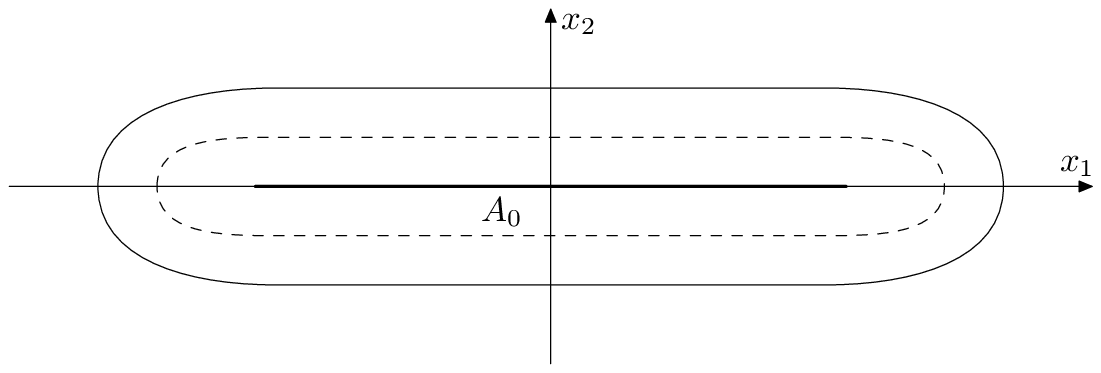}}   
\botcaption{Figure 5. two-step soul}
\endcaption

We observe that
$$l_0=\max\{d(p, \p\Omega_0)\,|\,p\in\Omega_0\}=1.$$
Let $\Omega_s=\{q\in\Omega_0\,|\,d(q,\p\Omega_0)\geq s\}$ and $r(q)=d(q, \p\Omega_0)$. The
maximum subset $A_0=r^{-1}(1)$ is an interval, i.e., $A_0=\{(x_1,0)\,|\,|x_1|\leq 4 \}$.
The soul of $\Omega_0$ is a single point set $A_1=\{(0,0)\}$. The level set $\p\Omega_s$
has strict convex points in $(\p\Omega_s)\cap[\R\times\{0\}]$ for $0\leq s<1$.

(b) We consider a family of rectangles $\{\Omega_s\}$, where
$$\Omega_s=\{(x_1,x_2)\,|\,|x_1|\leq 2-s, ~|x_2|\leq 1-s\}$$
for $0\leq s\leq 1$.

\vskip.3in

\epsfxsize=4.4cm             
\centerline{\epsfbox{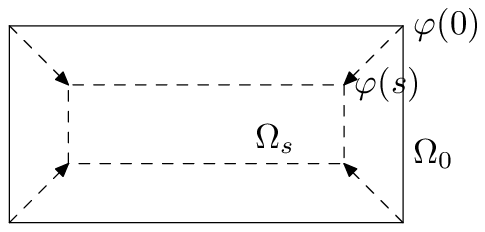}}   
\botcaption{Figure 6. Corner points and equi-distance evolution}
\endcaption

The corner point $q_s=(2-s,1-s)$ is a strictly convex point of
$\p\Omega_s$ for $0\leq s\leq 1$. It is clear that corner points
are preserved by Perelman-Sharafutdinov semi-flow of distance
functions.

(c) (Fermi coordinates and Riccati equation) Suppose that $\Omega_0$ is a convex domain
with a smooth boundary $\p\Omega_0$ in a smooth Riemannian manifold $M^n$ with non-negative
curvature. Suppose that $\varphi: [0, l]\to \Omega_0$ is a length-minimizing geodesic segment
of unit speed from $\p\Omega_0$ with
$$d(\varphi(s), \p\Omega_0) =s   \tag 4.2$$
for $0\leq s\leq l$.

\vskip.3in

\epsfxsize=5cm             
\centerline{\epsfbox{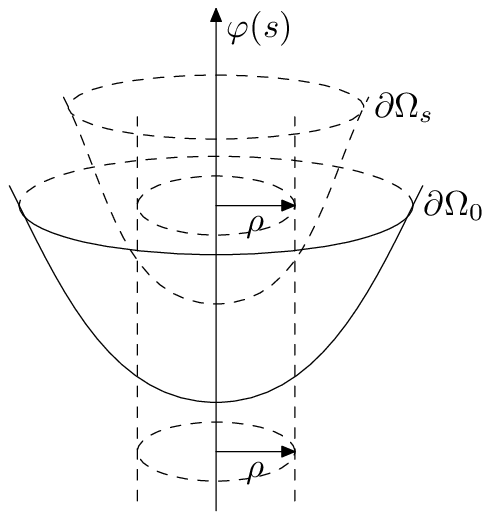}}   
\botcaption{Figure 7. fermi coordinates around a curve}
\endcaption

Let $\{\vec{E}_j(t)\}_{j=1}^n$ be a parallel orthonormal frame along the geodesic segment $\varphi$
such that $\vec{E}_n(t)=\varphi'(t)$. If $\p\Omega_t$is convex and if $\varphi'(t)$ is an inward
unit normal vector, then $\p\Omega_0$ lies at one side of $[M^n-H^{n-1}_{\varphi(t)}]$ where
$$H^{n-1}_{\varphi(t)}=\{\Exp_{\varphi(t)}\big(\rho\vec{E}(t)\big)\,|\,\vec{E}(t)\perp\varphi'(t)\}.$$
We consider the corresponding second fundamental form
$$\II_{\p\Omega_t}(X, Y)=-\li\varphi'(t), \nabla_XY\ri$$
for $\{X, Y\}\in T_{\varphi(t)}(\p\Omega_t)$. If
$$u_{ij}(t)=\II_{\p\Omega_t}(\vec{E}_i(t), \vec{E}_j(t))$$
and if $u(t)=\big(u_{ij}(t)\big)_{(n-1)\times(n-1)}$ is a matrix-valued representation of $\II_{\p\Omega_t}$
then the Riccati equation
$$u'+u^2+R=0   \tag4.2$$
holds, where
$$R_{ij}(t)=\li R(\varphi'(t), \vec{E}_i(t))\varphi'(t), \vec{E}_j(t)\ri,$$
$R(t)=\big(R_{ij}(t)\big)_{(n-1)\times(n-1)}$ is the curvature matrix function and
$$R(X, Y)Z = -\nabla_X\nabla_Y Z + \nabla_Y\nabla_X Z +\nabla_{[X,Y]} Z.$$
When $M^n$ has non-negative sectional curvature, the matrix $R(t)$ is positive semi-definitive:
$$R(t)\geq 0.  \tag 4.4 $$
It follows that
$$u'(t) = - u^2(t) - R(t) \leq 0   \tag 4.5$$
and hence the second fundamental form
$$\II_{\p\Omega_t}(E(t), E(t)) = \li u(t)E(t), ~E(t)\ri$$
is a monotone function of $t$ for any parallel vector field $\{E(t)\}$ along $\varphi$. This gives
rise to a version of Theorem $1.9$ for this special case. \qed

If $\varphi: [a, b]\to X^n$  is a Perelman-Sharafutdinov curve of
a distance function (not necessarily a length-minimizing geodesic)
in a singular space $X^n$ of non-negative curvature, there is no
Riccati equation available, we will use Trapezoid Comparison
Theorem (cf. Theorem $1.11$ above) instead.

\head \S 4.1. Angular estimates under equi-distance evolution in dimension $2$
\endhead

In this sub-section, we provide a proof of Theorem $1.9$ for the case when $X^n$ has dimension
$2$. We also derive some preliminary results for all dimensions.

Let us first recall the definition of barrier functions and barrier hypersurfaces, so that
we can estimates the concavity of Busemann function $h$ and its level sets $\p\Omega_c=h^{-1}(c)$.

\proclaim{Definition 4.2} \rom(Calabi \cite{Ca57} \rom) Let $h: X^n\to \R^n$ be a continuous
function. We say that $\Hess(h)(p)\leq \lambda$ if for any quasi-geodesic $\sigma: (-r,r)\to X^n$
with $\sigma(0)=p$ there exists an upper barrier function $\hat{h}$ such that $(\hat h \circ \sigma)(t)$
is smooth in $t$,
$$h(\sigma(0))=\hat{h}(\sigma(0)), ~~h(\sigma(t))\leq \hat{h}(\sigma(t))$$
and
$$\frac{d^2[\hat{h}(\sigma(t))]}{dt^2}\Big|_{t=0}\leq\lambda$$
hold.
\endproclaim

We remark that, in above definition, we need to choose two-sided barrier functions instead of
one-sided barrier functions to estimate the second derivative.

For example, if $f(t)$ is a smooth function of $t$, then we have the Taylor expansion
$$f(t)=f(0) +f'(0)t +\frac{1}{2}f''(0)t^2 +o(t^2)$$
where $\lim\limits_{t\to 0}\frac{o(t^2)}{t^2}=0$. It follows that
$$\align
f''(0) & = \lim_{t\to 0} \frac{f(t)+f(-t)-2f(0)}{t^2} \\
       & = \lim_{t\to 0} \frac{\frac{f(t)-f(0)}{t} -\frac{f(0)-f(-t)}{t}}{t}.
\endalign
$$
The second derivative $f''(0)$ measures the rate of change for slopes
$$\big[ \frac{f(t)-f(0)}{t} -\frac{f(0)-f(-t)}{t} \big] \sim \lambda t$$
up to the first order.

The geodesic curvature of a curve $\psi: [-r, r]\to \p\Omega_0$ in a smooth hypersurface
$\p\Omega_0$ of a smooth Riemannian manifold is given by
$$\II_{\p\Omega_0}(\psi'(0), \psi'(0)) = -\li\nabla_{\psi'}\psi', \varphi'(0)\ri
  =-\li\psi', \nabla_{\psi'}\vec{N}\ri\big|_{t=0}$$
where $\vec{N}$ is a smooth unit normal vector of $\p\Omega_0$.
The geometric quantity $-\li\psi', \nabla_{\psi'}\vec{N}\ri$
measure changes of slope as well.

We consider another example to indicate why we need to choose two sided barrier function.
Let us consider the function
$$h: \R^2 \to \R^2, ~~(x_1, x_2)\mapsto \frac{1}{\sqrt{2}}(x_2-|x_1|).$$

\vskip.3in

\epsfxsize=6cm             
\centerline{\epsfbox{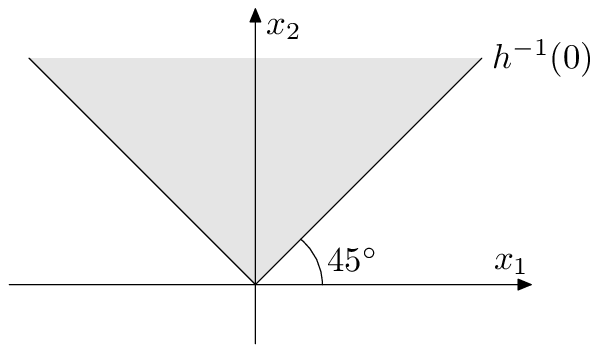}}   
\botcaption{Figure 8. non-smooth Busemann function}
\endcaption

Clearly, $h$ is a piecewise linear function. Hence, $h(x_1,x_2)$ has $\Hess(h)|_{p}=0$ if
$p=(x_1, x_2)$ satisfies $x_1\neq 0$. However, at $p=(0, 0)$, one can verify that
$$\Hess ( h )|_{(0, 0)} = -\infty$$
in barrier sense.

By our definition, $p=(0, 0)$ is a strictly convex point of $h^{-1}(0)$, where $\Omega_0=h^{-1}([0, +\infty))$.
In this example, $x_1$-axis lies below the level curve $h^{-1}(0)$. We need to derive some
preliminary results, in order to construct ``supporting cones" of $h^{-1}(c)$ outside a metric ball $B_r(p)$.

\proclaim{Proposition 4.3} \rom(Compare with \cite{CG72} \rom) Let $X^n$, $h(x)$ and $\{\Omega_s\}$
be as above. Suppose that $\sigma: (-\infty, +\infty) \to X^n$ is a quasi-geodesic of unit
speed with $\{\sigma(a), \sigma(b)\}\subset\p\Omega_s=h^{-1}(s)$ and $a<b$. Then the following hold.

$(4.3.a)$ If $\max\limits_{a\leq t\leq b} \{h(\sigma(t))\}> s$, then $h(\sigma(t))<s$ for any $t\notin[a, b]$;

$(4.3.b)$ If $\max\limits_{a\leq t\leq b} \{h(\sigma(t))\}= s$, then $h(\sigma(t))\leq s$ for all $t$.

Consequently, if $\tilde{\sigma}: [0, +\infty) \to \Omega_c$ is a quasi-geodesic segment with
$\{\tilde{\sigma}(0), \tilde{\sigma}(r)\}\subset\p\Omega_c$ and $r >0$, then
$$\tilde{\sigma}([r, +\infty)) \subset [X^n -\Int(\Omega_c)] =h^{-1}((-\infty, c]).$$
\endproclaim

\demo{Proof} Let $\eta(t)=h(\sigma(t))$. by our assumption, we
know that $\eta(t)$ is a concave function. Our conclusion are
direct consequences of concavity of $\eta(t)$. \qed
\enddemo

\vskip.3in

\epsfxsize=5cm             
\centerline{\epsfbox{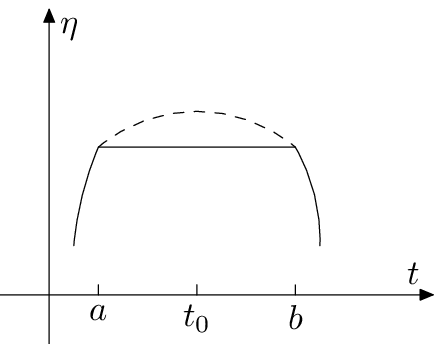}}   
\botcaption{Figure 9}
\endcaption

On a singular space $X$ with curvature $\geq -1$, the extension of quasi-geodesic segment
$\sigma: [a, b] \to X^n$ to a ``longer" quasi-geodesic $\tilde{\sigma}: \R\to X^n$ is not
necessarily unique. Thus, it is necessary for us to recall the notion of the gradient exponential
map (cf. \cite{Petr07}).

Such a gradient exponential map $\Exp_p : T_pX \to X$ is related to the gradient flow of the
distance function $r(x)=d(p, x)$ which we now describe.

If $X$ is an Alexandrov space with metric $d$,then we denote by
$\lambda X$ the space $(X, \lambda d)$. Let $i_{\lambda}: \lambda
X \to X$ be the canonical map. The Gromov-Hausdorff limit of
pointed spaces $\{(\lambda X, x)\}$ as $\lambda\to+\infty$ is the
tangent cone of $(T_x(X), o_x)$ of $X$ at the point $x$, (see \S
7.8.1 of \cite{BGP92}. For any $\lambda$-concave function, the
function $d_x f: T_x(X)\to\R$ defined by
$$d_x f = \lim_{\lambda\to+\infty} \frac{f\circ i_{\lambda} -f(x)}{1/\lambda}$$
is called the differential of $f$ at $x$.

The gradient vector $\nabla f$ of $f$ at $x$ is related to the inner product of $T_x(X)$. For
any pair of vectors $\vec{u}$ and $\vec{v}$ in $T_x(X)$, we define
$$\li\vec{u}, \vec{v}\ri = \frac{1}{2}(|\vec{u}|^2 +|\vec{v}|^2 -|\vec{u}\vec{v}|^2) = |\vec{u}||\vec{v}|\cos\theta$$
where $\theta$ is the angle between $\vec{u}$ and $\vec{v}$, $|\vec{u}\vec{v}|=d_{T_x(X)}(\vec{u}, \vec{v})$,
$|\vec{u}|=d_{T_x(X)}(\vec{u}, o)$, and $o$ denotes the origin of the tangent cone. It follows that
$$\cos\theta = \frac{\li\vec{u}, \vec{v}\ri}{|\vec{u}||\vec{v}|}.$$

\proclaim{Definition 4.4} \rom(\cite{Per91}, \cite{petr07} \rom) For any given semi-concave function
$f$ on $X$, a vector $\vec{\eta}\in T_x(X)$ is called a gradient of $f$ at $x$ ($\eta=\nabla f$ in
short) if the following holds:

$(i)$ $d_x f(\vec{v}) \leq \li\vec{\eta}, \vec{v}\ri$ for any $\vec{v}\in T_x(X)$;

$(ii)$ $d_x f(\vec{\eta}) = |\vec{\eta}|^2$.

\endproclaim

It is known that any semi-concave function has a uniquely defined gradient vector field.
Moreover, if $d_x f(\vec{v})\leq 0$ for all $\vec{v}\in T_x(X)$, then $\nabla f|_x=0$. In
this case, $x$ is called a critical point of $f$. Otherwise, we set
$$\nabla f = (d_x f)(\vec{\xi})\vec{\xi}$$
where $\vec{\xi}$ is the (necessarily unique) unit vector for which $d_x f$ attains its
positive maximum on $\Sigma_x(X)$ and $\Sigma_x(X)$ is the space of directions of $X$ at
$x$.

When $X^n$ has curvature $\geq 0$, its energy function $f(x)=\frac{1}{2}[d(x, q)]^2$ is
$1$-concave. We consider the semi-flow
$$\Phi_f^t(q) = \alpha_q(t)$$
where $\alpha_q(0)=q$ and
$$ \frac{d^+ \a_q(t)}{dt} =\nabla f\big|_{\a_q(t)}.$$
Recall that
$$\lim_{\l\to +\infty} (\l X, q) = (T_q(X), o_q).$$
We define $\gexp_q: T_q(X)\to X$ by
$$\gexp_q =\lim\limits_{t\to +\infty} \Phi_f^t \circ i_{e^t},$$
where $i_\l: \l X \to X$ is the canonical map. In fact, for each unit direction
$\vec\xi\in \Sigma_p(X)$, the radial curve $\a_{\vec\xi}: t\mapsto \gexp_p (t\vec\xi)$ satisfies
the equation
$$ \frac{d^+ \a_{\vec\xi}}{dt} =\frac{\hat r_p(\a_{\vec\xi}(t))}{t} \nabla \hat r_p$$
where $\hat{r}_p(x)=d(x, p)$.

\proclaim{Theorem 4.5} \rom(\cite{PP94}, \cite{Petr07} p.152\rom) Let $X$, $\hat r_p$ and $\gexp_p$
be as above. If $\vec \xi\in T_p (X)$ and if $h$ is a concave function, then
$$ h(\gexp_p(t\vec\xi)) \leq h(p) +t(d_p h)(\vec\xi)$$
for all $t\geq 0$.
\endproclaim

In our Theorem $4.5$, there is a first order term $(d_p
h)(\vec\xi)$. Thus, we need to recall the first variational
formula.

\proclaim{Theorem 4.6} \rom(\cite{BBI01}, p.125\rom) Let $X^n$ be a complete Alexandrov space of
non-negative curvature and $f_A(x)=d(x, A)$, where $A$ is closed subset of $X$. Suppose that
$x\notin A$, $\Uparrow_x^A$ is the set of minimal directions from $x$ to $A$, $|\vec{u}|=1$ and
$$\theta=\min\{\ang_x(\vec{u}, \vec{w})\,|\,\vec{w}\in\Uparrow_x^A\}.$$
Then $df_A(\vec{u})=-\cos\theta$.
\endproclaim

Among other things, we need to use Theorem $1.11$ to verify Theorem $1.9$.

\demo{Proof of Trapezoid Comparison Theorem $1.11$} Let us make some observations on Figure 3.
When curvature is non-negative, the (quasi)-geodesic triangle with vertices $\{ \hat q, p_2, p_3\}$
satisfies
$$
\big[d(p_2, \varphi_3(l\cos\alpha_2))\big]^2 \leq [\ell^2 + (\ell \cos \alpha_2  )^2 - 2 \ell^2 (\cos \alpha_2  )^2  ] = (l\sin\alpha_2)^2.
$$
In order to complete our proof, it is sufficient to estimate both $l$ and $\alpha_2$ from above.

Since $X^n$ has non-negative curvature, by the triangle comparison theorem we have
$$\align
l^2 & =[d(\hat{q}, p_2)]^2 \leq \big(\frac{s}{\sin\hat{\beta}}\big)^2
       + r^2 -2\frac{rs}{\sin\hat{\beta}}\cos(\pi -\hat{\beta}) \\
    & = \frac{s^2}{(\sin\hat{\beta})^2} + r^2 + 2rs\cot\hat{\beta} \\
    & = \hat{l}^2.   \tag 4.6
\endalign
$$

Using Euclidean trigonometry, for the comparison triangle $\D_{\hat{q}\hat{p}_2\hat{p}_3}$ with
vertices $\{\hat{q}, \hat{p}_2, \hat{p}_3\}$ in $\R^2$ (see Figure 10 below), we have the inequality
$$\align
\big[d(\hat{p}_2, \hat{\varphi}_3(\hat{l}\cos\hat{\alpha}_2))\big]^2
& = \hat{l}^2 + (\hat{l}\cos\hat{\alpha}_2)^2 -2\hat{l}^2(\cos\hat{\alpha}_2)^2 \\
& = (\hat{l}\sin\hat{\alpha}_2)^2 \leq s^2,  \tag 4.7
\endalign
$$
where
$$\hat{\alpha}_2 = \hat{\beta} -\hat{\alpha}_1.  \tag 4.8$$

\vskip.3in

\epsfxsize=6.2cm             
\centerline{\epsfbox{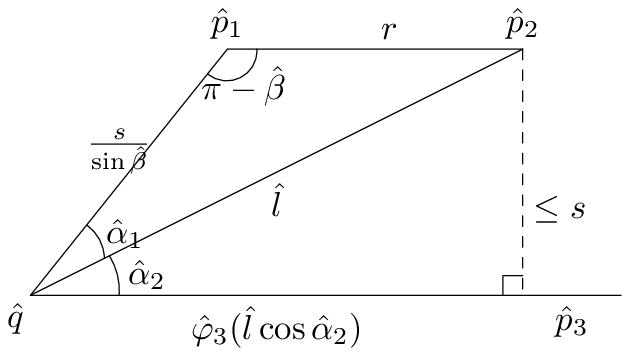}}   
\botcaption{Figure 10}
\endcaption

Let us now consider the Euclidean triangle $\D_{\hat{q}\hat{p}_1\hat{p}_2}$ of the side
lengths $\{\hat{l}, r, \frac{s}{\sin\hat{\beta}}\}$. Since $X^n$ has non-negative curvature,
$\D_{\hat{q}p_1p_2}$ is fatter than $\D_{\hat{q}\hat{p}_1\hat{p}_2}$. It follows that
$$\alpha_1\geq \hat{\alpha}_1 ~~\text{and} ~~\alpha_2=\hat{\beta}-\alpha_1\leq \hat{\beta}-\hat{\alpha}_1=\hat{\alpha}_2.  \tag4.9$$
Therefore, by $(4.6)-(4.9)$ we conclude that $l \le \hat l$, $\alpha_2 \le \hat{\alpha}_2$ and
$$\big[d(p_2, \varphi_3(l\cos\alpha_2))\big]^2 \leq (l\sin\alpha_2)^2
  \leq (\hat{l}\sin\hat{\alpha}_2)^2 \le s^2.$$
This completes the proof of Theorem $1.11$.
\qed\enddemo

We conclude this subsection by a special case of Theorem $1.9$ when $\dim(X^n)=2$.

\proclaim{Theorem 4.7} Let $X^2$ be an open Alexandrov surface with non-negative curvature,
$$h(x)=\lim\limits_{t\to +\infty} [d(x, \p B_t(x_0))-t]$$
and $\W_s =h^{-1}([s,+\infty))$. Suppose that $\varphi: [0,l] \to X^2$ is a Perelman-Sharafutdinov
curve for $h$ defined by
$$ \frac{d^+ \varphi}{dt} =\frac{\nabla h}{|\nabla h|^2} \Big|_\varphi,$$
and $\eta(s,r) =\theta^{\Omega_{h(\varphi(s))}}_{\varphi(s),\vec{u}(s)}(r)$  where $\vec{u}(s)=\frac{d^{-}\varphi}{ds}$
is the normalized left derivative of $\varphi$. Then
$$ \frac{\p\eta}{\p s} (s,r) \geq 0.  \tag4.10$$
\endproclaim

\demo{Proof} If $\W_c$ is a totally convex subset of $X^2$ and if $\dim (\W_c)=2$, then Perelman
(cf. \cite{Per91}, Chapter 6) calculated $\nabla h$ as follows.

\vskip.3in

\epsfxsize=6cm            
\centerline{\epsfbox{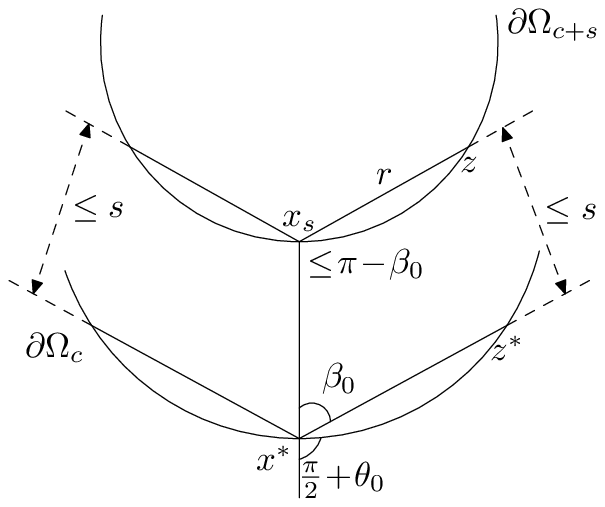}}   
\botcaption{Figure 11}
\endcaption
\ms

If $x\in\p\W_c$ and $\beta =\frac12 \roman{diam}(\Sigma_x(\W_c))$, then Perelman
(cf. \cite{Per91} page 33, line 1) observed that
$$ |\nabla h(x)| =\sin [\beta(x)]. \tag 4.11 $$
For fixed $r$, we would like to show  $\eta(s,r) =\o^{\W_{h(\varphi(s))}}_{\varphi(s)}(r)$
is an increasing function in $s$.

Perelman (cf. \cite{Per91}, Chapter 6) already proved that the Sharafutdinov retraction is
distance non-increasing. Hence the map
$$ \pi: \p\W_c \to \p \W_{c+\e}$$
is a distance non-increasing map.

If $x, y\in \p \W_{c+\e}$,  $x^* \in \pi^{-1}(x)$ and $y^* \in\pi^{-1}(y)$, then
$$  d(x^*,y^*) \geq d(x,y) . $$

Since $X^2$ has non-negative curvature, it is known that diam$[\Sigma_x^1(X^2)] \leq \pi$.

If $x^*=\varphi(0)$ and $\o_0(r)=\o^{\W_{h(x^*)}}_{x^*} (r)$, then
$$ \beta_0 (r) \leq \pi -(\frac\pi{2} +\o_0) =\frac\pi{2} -\o_0. \tag 4.12$$
Let $x_s=\varphi(s)$. We consider the left derivative $\vec v_s$ of $\varphi$ at $x_s$.
Since $x_s$ is the nearest point on $\p\W_{h(x_s)}$ to $x_0$, the vector $\vec v_s$ is at
least normal to $\p \W_{h(x_s)}$ at $x_s$. By our assumption, we have
$$ \ang_{x^*} (x_s, [\p\W_{h(x^*)} -B_r(x^*)]) \leq \frac\pi{2}
-\o^{\W_{h(x^*)}}_{x^*} (r) +o(s), \tag 4.13 $$
where $\lim\limits_{s\to 0} \frac{o(s)}{s} =0$.

We now consider any given quasi-geodesic $\sigma_0: [0,+\infty) \to X^2$ such that
$\sigma_0(0) =x^* \in \p\W_c$ and $ \sigma_0(l)\in \p\W_c$ with $l>0$,
$$ d(x^*, \sigma_0(l)) \leq r. \tag 4.14$$
It follows from Proposition $4.3$ that
$$ \sigma_0([l,+\infty)) \subset h^{-1}((-\infty, c^*]) =[X^2 -\roman{int}(\W_{c^*})]. $$
By our discussion above, we have
$$ \beta^*_0 =\ang_{x^*} (x_s, \sigma'_{0}(0)) \leq \frac\pi{2} - \o_{x^*}^{\W_{c^*}} (r) +o(s). \tag 4.15$$

Our technical goal is the following

\ni {\bf Claim 4.7a}. {\it  Let $X^2$, $h$, $\{ \W_c\}$, $\o_x^\W (x)$,
$\sigma_0$ and $\beta_0^*$ be as above (see Figure 11). Suppose that $0 < s \le \frac{r}{32}$
and that $\sigma_s: [0,+\infty) \to X^2$ is a quasi-geodesic with $\sigma_s(0)=x_s$,
$\sigma_s(l) \in [\W_{c^*}-B_r(x_s)] \cap U_{\frac{r}{8}} (\sigma_0([l,+\infty)))$ and
$$ \a=\ang_{x_s}(x^*, \sigma'_s(0)) \leq \pi -\beta_0^* =\frac\pi{2} +\o_s^*. \tag 4.16$$
Then
$$ \sigma_s([l, +\infty)) \subset U_s (\sigma_0([l,+\infty))) \subset h^{-1}
((-\infty, c^* +s]). \tag 4.17$$}

\ms If Claim $4.7a$ is true, then by (4.16)-(4.17), whenever $y \in [\W_{c^* +s} - B_r( x_s)]  \subset h^{-1}
((-\infty, c^* +s]) $,  we must have
$$ \o_{x_s}^{\W_{c^* +s}}(r) \geq \pi-\beta^*_0 =\frac\pi{2} +\o_{x^*} +o(s), \tag 4.18$$
where $ \lim_{s \to 0}\frac{o(s)}{s} = 0$.
It follows that if $\eta(s,r) =\o_{x_s}^{\W_{c^* +s}}(r)$, then
$$ \frac{\p \eta(s,r)}{\p s} \geq 0 \tag 4.19$$
and Theorem $4.7$ holds.

It remains to verify Claim $4.7a$. Inspired by Perelman's formula,
$|\nabla h|(x^*) \sim \sin (\hat\beta_{x^*})$, where
$$ \hat \beta_{x^*} =\max\{ \ang_{x^*} (\vec w, T_x(\p\W_{c^*})) \}. $$
For a fixed $\beta_0(r)$, we choose a geodesic segment $\psi: [0,
l_s] \to \W_{c^*}$ of unit speed with $\psi(0) =x^*$, $\ang_{x^*}
(\sigma'_0(0), \psi'(0)) =\beta_0(r)$, $l_s=\frac{s}{\sin\beta_0}$ and $x_s = \psi(l_s)$.

Recall that $X^2$ has non-negative curvature, using triangle comparison theorem, we have
$$ \align
[d(x_s, \sigma_0(s\cot \beta_0))]^2 &\leq \frac{s^2}{(\sin\beta_0)^2}
+s^2 (\cot\beta_0)^2 -2s^2 \cot\beta_0 \frac1{\sin\beta_0} \cos\beta_0 \\
&=s^2 [(\frac1{\sin\beta_0})^2 -(\cot\beta_0)^2] \\
&=s^2.
\endalign
$$
It follows that
$$ d(x_s, \sigma_0(\R)) \leq s. $$
Let us choose $y^* \in\sigma_0(\R)$ such that
$$ d(x_s, y^*)=d(x_s, \sigma_0(\R)). \tag 4.20$$

Since $X^2$ has non-negative curvature, the sum of interior angles of any geodesic triangle
is greater than or equal to $\pi$. It follows that
$$ \ang_{x_s} (x^*, y^*) \geq \pi-(\frac\pi{2}+\beta_0) =\frac\pi{2} -\beta_0.$$
For each $z\in [\W_{c^*+s}-B_r(x_s)] \cap U_{\frac{r}{8}} (\sigma_0([l_s,+\infty)))$, there
are two cases.

If $\ang_{x_s} (z,y^*) \geq \frac\pi{2}$, then
$$ \ang_{x_s} (z,x^*) \geq \pi -\beta_0 =\frac\pi{2} +\o^*,$$ and we are done.

If  $\ang_{x_s} (z,y^*) < \frac\pi{2}$, we get a contradiction as follows. Let
$z^* \in \sigma_0((0,+\infty))$ satisfy $d(z,z^*)=d(z, \sigma_0((0,+\infty)))$. Since the
Sharafutdinov flow is a distance non-increasing map, we must have
$$ d(x^*,z^*) \geq d(x_s,z) \geq r.$$
By Proposition $4.3$, we have $z^* \not\in \roman{Int}(\W_{c^*})$. It follows that
$$ h(z^*) \leq c^*. \tag 4.21$$
When $\ang_{x_s} (z,y^*) < \frac\pi{2}$, by Theorem $1.11$, we
would have
$$ d(z,z^*) <s. \tag 4.22$$
This together (4.21) would imply that
$$ h(z) <c^* +s \tag 4.23$$
which contradicts to $z\in \p\W_{c^*+s}$. This completes the proof
of Claim $4.7a$ and hence Theorem $4.7$. \qed
\enddemo

\head \S 4.2. Parallel supporting cones outside a $\rho$-tube and generalized Fermi coordinates
\endhead

In this subsection, we elaborate the proof of Theorem $4.7$ and extend it to the higher dimensional
case with necessary modifications.

\vskip.3in

\epsfxsize=6.2cm             
\centerline{\epsfbox{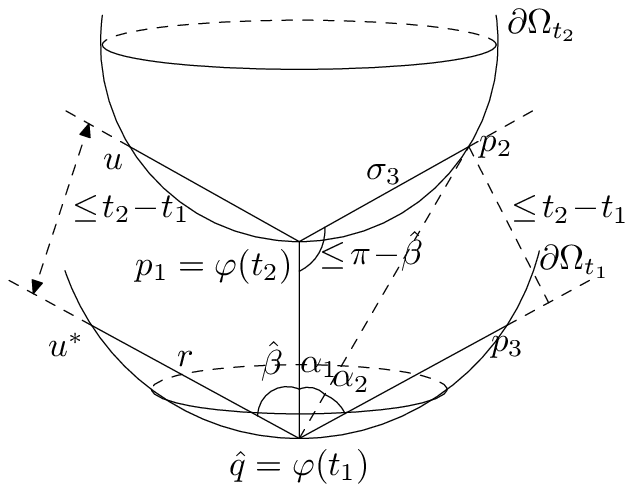}}   
\botcaption{Figure 12. Parallel supporting cones outside $r$-tubes}
\endcaption

Let us first make sure hat for each $p\in X^n$, there exists $r_0$ such that any quasi-geodesic
segments $\sigma: [0, +\infty) \to X^n$ starting from $p$ will leave $B_{r_0}(p)$, i.e.,
$\sigma(t)\notin B_{r_0}(p)$ for relative large $t$.

\proclaim{Proposition 4.8} \rom(Perelman \cite{Per94b} \rom) Let $X^n$ be a complete Alexandrov
space with curvature $\geq -1$. For each $\hat{x}\in X^n$, there exist positive numbers $\lambda$,
$\delta$ and a strictly concave function $f: B_{\delta}(p) \to (-\infty, 0]$ such that

$(i)$ $f(p)=0$;

$(ii)$ $B_{\frac{\varepsilon}{\lambda}}(p)\subset f^{-1}([-\varepsilon, 0])\subset B_{\lambda\varepsilon}(p)$
for $\varepsilon<\frac{\delta}{4\lambda}$, where $\lambda$ and $\delta$ depend on $p$.

\endproclaim

\proclaim{Corollary 4.9} Let $X^n$ be a complete Alexandrov space with curvature $\geq -1$.
Suppose that $p\in X^n$, $\lambda$, $\delta$ and $f$ be as in Proposition $4.8$ above. Then
any quasi-geodesic $\psi: [0, +\infty)\to X^n$ with $\psi(0)=p$ must leave $B_{\frac{\delta}{\lambda}}(p)$
when $t$ is sufficiently large.
\endproclaim

\demo{Proof}  Since $f$ is strictly concave on $B_{\delta}(p)$ and $p$ is a unique maximum point
of $f$, there exists $\e_0>0$ such that
$$-\e_1 = f(\psi(\e_0)) < f(\psi(0)) =0.$$
It follows from the concavity of $f$ that for $t>\e_0>0$, we have
$$\frac{f(\psi(t)) -f(\psi(\e_0))}{t-\e_0} \leq \frac{f(\psi(\e_0)) -0}{\e_0 -0}$$
and hence
$$f(\psi(t)) \leq \frac{f(\psi(\e_0))}{\e_0}(t-\e_0) = -\frac{\e_1}{\e_0}(t-\e_0).$$
Therefore, for sufficiently large $t$, we have
$$\psi(t) \notin f^{-1}([-\e_0, 0]).$$
Recall that $B_{\frac{\e_0}{\lambda}}(p)\subset f^{-1}([-\e_0, 0])$. Hence we have
$\psi(t)\notin B_{\frac{\e_0}{\lambda}}(p)$ for sufficiently large $t$.
\qed
\enddemo

We now discuss the supporting hypersurface when $r=0$ and $\theta_{p,\vec{u}}^{\Omega_c}(0)=\frac{\pi}{2}$.

\proclaim{Theorem 4.10}\rom(\cite{Per91} Chapter 6, \cite{Petr07} p.156, \cite{CDM09} p.40\rom)
Let $X^n$ be an open complete Alexandrov space with non-negative curvature,
$$h(x)=\lim\limits_{t\to +\infty} [d(x,\p B_t(x_0))-t]$$
and $\W_c=h^{-1}([c,+\infty))$. Suppose that $\hat x\in \W_{\hat c}$,
$\vec v\in \Uparrow_{\hat x}^{\p\W_{\hat c-\delta}}$, and $\e$ is given by Proposition $4.8$.
We consider $ H^{n-1}_{\hat x,\e} =\{ y\in B_\e(\hat x) | \sigma:[0,l] \to B_\e (\hat x) \ \text{is a quasi-geodesic},\
\sigma(0)=\hat x,\ \sigma'(0) \perp \vec v,\  l<\e,\  y=\sigma(l) \}$.
Then $\p\W_{\hat c} \cap B_\e (\hat x)$ lies above $H_{\hat x,\e}^{n-1}$ \rom(i.e., $h(y)\leq \hat c$
for $y \in H^{n-1}_{\hat x,\e}$.\rom)
\endproclaim

\demo{Proof} It is known (cf. \cite{Wu79})that
$h(y)=d(y,\p\W_{\hat c-\delta}) +\hat c$ for $y\in \W_{\hat c-\delta}$.

\vskip.3in

\epsfxsize=5.5cm             
\centerline{\epsfbox{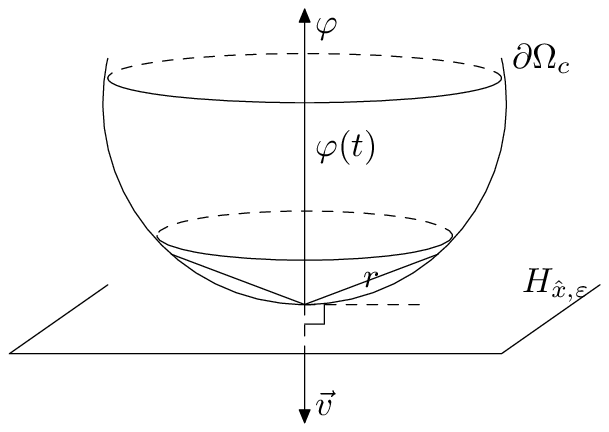}}   
\botcaption{Figure 13. Supporting hypersurfaces with $r=0$}
\endcaption
\ms

Perelman also showed that if $N=\p\W_{\hat c-\delta}$, then $r_N(y)=d(y,N)$ is concave for
$y\in \W_{\hat c -\delta}$.  For any quasi-geodesic $\psi: [0,+\infty) \to X^n$ with
$\psi(0)=\hat x$, $\psi'(0)=\vec w$, and
$$ \ang_x(\vec v,\vec w) \leq \frac\pi{2},\tag 4.24$$
we have
$$ \frac{d[h(\psi(t))]}{dt} \Big|_{t=0}  \leq 0. $$
By an equivalent definition of quasi-geodesics, $t\mapsto h(\psi(t))$ is a concave function.
Therefore, we have
$$ h(\psi(t)) \leq h(\psi(0)) =\hat c. \qed$$
\enddemo

We now consider supporting cones of $\Omega_c$ outside the metric ball $B_r(p)$ related to
angular excess function $\theta_{p,\vec{u}}^{\Omega_c}(r)>0$ for $r>0$.

\vskip.3in

\epsfxsize=6.2cm             
\centerline{\epsfbox{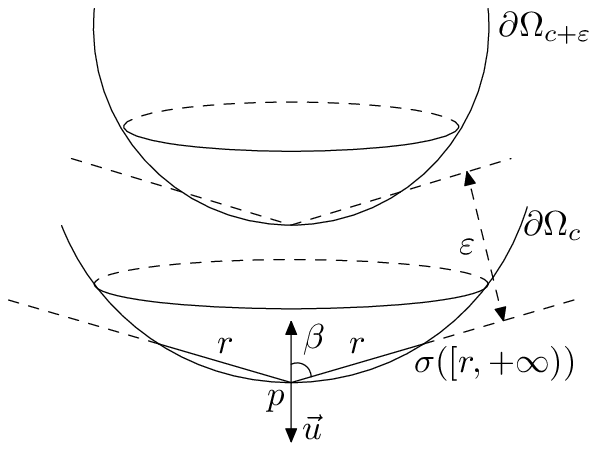}}   
\botcaption{Figure 14. Supporting cones outside the ball $B_r(p)$ }
\endcaption

In a smooth Riemannian manifold $M^n$, its tangent cone $T_p(M^n)$ is isometric to the Euclidean
space $\R^n$. Hence, any unit vector $\vec{u}$ has its anti-podal vector $-\vec{u}$. For
singular space $X^n$ with curvature $\geq -1$, we need to recall the notion of polar vectors, so
that we can discuss relations between inner angles and outer angles.

\proclaim{Definition 4.11} $(i)$ Two vectors $\vec{u}$, $\vec{v}\in T_p(X)$ \rom(not necessarily
the same length\rom) are called polar if for any vector $\vec{w}\in T_p(X)$, the inequality
$$\li \vec{u}, \vec{w}\ri + \li\vec{v}, \vec{w}\ri \geq 0$$
holds.

$(ii)$ A vector $\vec{w}\in T_p(X^n)$ is called a supporting vector of a semi-concave function
$f$ at $p$ if
$$(d_pf(\vec{\xi})) \leq -\li \vec{\xi}, \vec{w}\ri$$
holds for any $\vec{\xi}\in T_p(X^n)$.
\endproclaim

It is known that the set of supporting vectors of $f$ at $p$ form a non-empty convex subset of
$T_p(X^n)$, see \cite{Petr07} p143. In fact, if $r(p)=d(p, A)$ with $p\notin A$ and if
$\vec{w}\in\Uparrow_p^A$, then $\vec{w}$ is a supporting vector of $r_A$ at $p$, by the first
variational formula.

Using the existence of quasi-geodesics, Petrunin (cf. \cite{Petr07} p194) showed that for unit
vector $\vec{u}\in\Sigma_p(X^n)$, there exists a polar unit vector $\vec{u}^{*}\in\Sigma_p(X^n)$
such that
$$\li \vec{u}, \vec{w}\ri + \li\vec{u}^{*}, \vec{w}\ri \geq 0$$
for any $\vec{w}\in T_p(X^n)$. Consequently
$$\ang_p(\vec{u}, \vec{w}) + \ang_p(\vec{u}^{*}, \vec{w}) \leq \pi   \tag 4.25$$
for any $\vec{w}$.

We now ready to construct the supporting cone of $\Omega_c$ outside $B_r(p)$. Recall that
if $p\in\p\Omega_c$ and $\vec{u}$ is a unit vector at least normal to $\Omega_c$ at $p$ then
$$\theta_{p,\vec{u}}^{\Omega_c}(r) =\inf \{\ang_p(\vec{u}, x)|x\in[\Omega_c-B_r(p)]\} -\frac{\pi}{2}.$$

\proclaim{Proposition 4.12} Let $X^n$, $h$, $\Omega_c=h^{-1}([c, +\infty))$, $p\in\Omega_c$ and
$\theta_{p,\vec{u}}^{\Omega_c}(r)$ be as above. Suppose that $r>0$,
$\theta=\theta_{p,\vec{u}}^{\Omega_c}(r)>0$, $\beta=\frac{\pi}{2} -\theta$ and $\sigma: [0, l]\to X^n$
is geodesic segment with
$$\ang_p(\vec{u}^{*}, \sigma'(0)) \geq \beta   \tag4.26$$
where $\{\vec{u}, \vec{u}^{*}\}$ are polar vectors in $T_p(X^n)$. Then
$$\sigma([r, +\infty))\subset [X^n -\Int(\Omega_c)] = h^{-1}((-\infty, c]).   \tag4.27$$
Consequently, $\sigma([r, +\infty))$ lies below $\p\Omega_c$ relative to the vector $\vec{u}^{*}$.
\endproclaim

\demo{Proof}  This is a direct consequence of concavity for
Busemann function $h$. \qed \enddemo

\vskip.3in

\epsfxsize=6.2cm             
\centerline{\epsfbox{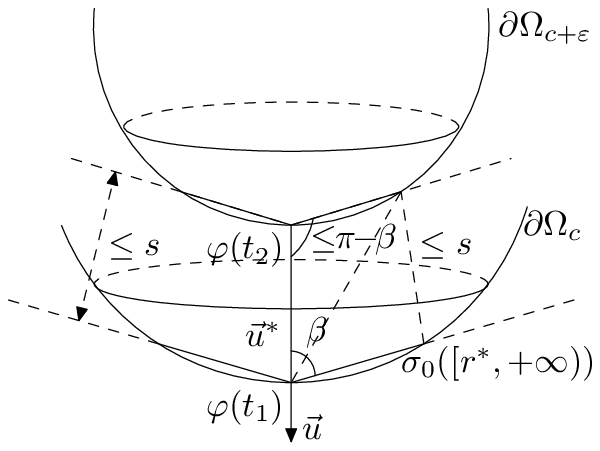}}   
\botcaption{Figure 15. Parallel cones along a curve}
\endcaption

\demo{Proof of Theorem $1.9$}  Let $\tilde{\varphi}: [0, l]\to X^n$ be a gradient exponential
radius curve given by
$$\tilde{\varphi}(s) =\gexp_p(\frac{s}{\sin\beta}u^{*}),$$
where $\beta=\frac{\pi}{2}-\theta_{p,\vec{u}}^{\Omega_c}(r) < \frac{\pi}{2}$.

Let us consider a conic hypersurface
$$\align
H_{\varphi(0),\beta}(r) = \{ \sigma_0(t)\,|\,& d(p, \sigma_0(t))\geq r, ~\sigma_0 ~\text{is a} \\
 & ~\text{quasi-geodesic with} ~\sigma_0(0)=p ~\text{and} ~\ang_p(\vec{u}^{*}, \sigma_0'(0))=\beta \}.
\endalign
$$
It follows from Proposition $4.12$ that the conic hypersurface $H_{\varphi(0), \beta}(r)$
lies below $\p\Omega_c$ relative direction $\vec{u}^{*}$. More precisely,
$$H_{\varphi(0), \beta}(r) \subset h^{-1}((-\infty, c]).  \tag 4.28$$
By our assumption $\tilde{\varphi}(s)=\gexp_p(\frac{s}{\sin\beta}\vec{u}^{*})$, we have
$$d(\tilde{\varphi}(s), \sigma_0(s\cot\beta)) \leq s + o(s)  \tag4.29$$
where $\lim\limits_{s\to 0}\frac{o(s)}{s}=0$.

We now consider a ``parallel" transport of $H_{\varphi(0), \vec{u}^{*}}(0)$ along the quasi-geodesic
$\varphi$ in the sense of Petrunin (cf. \cite{Petr98})
$$\align
\tilde{H}_{\tilde{\varphi}(s), \beta}(r) = \{x | & \ang_{\tilde{\varphi}(s)}(\varphi(0), x) = \pi -\beta,
~\exists \sigma_0 ~\text{such that} \\
  & \sigma_0(t)\in H_{\varphi(0), \beta}(r), \ang_{\tilde{\varphi}(0)}(\vec{u}^{*}, x)
    + \ang_{\varphi(0)}(x, \sigma_0'(0)) =\beta \}.
\endalign
$$

Recall that $h$ is a Busemann function
$$h(x) = (c-100) + d(x, \p\Omega_{c-100})$$
which is also distance function. Moreover, Perelman-Sharafutdinov semi-flow related to $h$ is
a distance non-increasing semi-flow. Using Theorem $1.11$ (Trapezoid Comparison Theorem), we see that,
for each $x\in \tilde{H}_{\varphi(s), \beta}$, we have
$$d(x, H_{\varphi(0), \beta}(r)) \leq s + o(s).   \tag4.30$$
If follows from $(4.28)$ and $(4.30)$ that
$$\tilde{H}_{\tilde{\varphi}(s), \beta}(r) \subset h((-\infty, c+s+o(s))).  \tag4.31$$
Thus, $\tilde{H}_{\tilde{\varphi}(s), \beta}(r)$ lies below
$\Omega_{c+s+o(s)}$ and outside $B_r(\tilde{\varphi}(s))$.

If $\vec{u}(s)$ is the left derivative of $\tilde{\varphi}$ at $s$, then
$$\align
\ang_{\tilde{\varphi}(s)}(\vec{u}(s), x) & \ang_{\tilde{\varphi}(s)}(\tilde{\varphi}(0), x) = \beta \\
  & \geq \pi -[\frac{\pi}{2} - \theta_{\tilde{\varphi}(0),\vec{u}}^{\Omega_c}(r)] \\
  & = \frac{\pi}{2} + \theta_{\tilde{\varphi}(0)}^{\Omega_c}(r)  \tag 4.32
\endalign
$$
Let $\theta_s(r)=\theta_{\varphi(s),\vec{u}(s)}^{\Omega_{h(\varphi(s))}}(r)$. Recall that if
$X^n$ has non-negative curvature, then $\gexp_q: T_q(X) \to X^n$ is a distance non-increasing
map. With extra efforts and using discussions above, we have
$$\theta_s(r) \geq \theta_0(r) + o(s)$$
where $\lim\limits_{s\to 0} \frac{o(s)}{s}=0$. It follows that
$$\frac{\partial\theta_s(r)}{\partial s}(0) \geq 0.  \qed$$
\enddemo

\head \S 5. Proof of Main theorem
\endhead

Let $a_0=\max\{h(x)\,|\,x\in X^n\}$ and $A_0=h^{-1}(a_0)$. Suppose that $X^n$ has positive
curvature on $B_{2\e_0}(\hat{x})$. Then using Proposition $1.8$, we see that
$$\theta_{\hat{x},\vec{u}}^{\Omega_c}(x) \geq ar >0   \tag 5.1$$
for some $a>0$ and all sufficiently small $r>0$. Suppose that $\varphi:[0,l]\to X^n$ is a
Perelman-Sharafutdinov curve for $h$ such that
$$
\cases
\frac{d^{+}\varphi}{dt}  = \frac{\nabla h}{|\nabla h|^2}, \\
\varphi(0)  = \hat{x}, \\
\varphi(l)  \in A_0
\endcases
$$
by Theorem $1.9$ and $(5.1)$, we see that
$$\theta_{\varphi(l), \vec{u}(l)}^{\Omega_{a_0}}(r) > ar >0   \tag 5.2$$
for all $0<r\leq r_0$. It follows from Proposition $1.8$ and $(5.2)$ that either $A_0$ is a
point set or $A_0$ has non-empty boundary with a ``strictly convex" boundary point $\varphi(0)\in\p A_0$.

When $\p A_0\neq\emptyset$, we let $\Omega_{a_0+s} =\{z\in A_0|d(z, \p A_0)\geq s\}$,
$l_i=\max\{d(z, \p A_0)|z\in A_0\}$ and $A_1=\Omega_{a_0+l_1}$. Using the same argument as above,
we can find $z_1\in\p A_1$ and $\vec{u}_1\in T_{z_1}(X)$ such that
$$\theta_{z_1,\vec{u}_1}^{A_1}(r) \geq a_1 r >0$$
for $r>0$ unless $\dim(A_1)=0$. Observe that
$$\dim(X^n) > \dim(A_0) > \dim(A_1)>\cdots.$$
Repeating above procedure $m\leq n$ times, we conclude that $A_m$ is of dimension zero. Therefore,
$X^n$ can be contractible to a point $z_m=A_m$ via multiple step Perelman-Sharafutdinov retraction.
\qed

Finally, we remark that our angular excess estimates are related to the length-excess function of
Alexander-Bishop (cf. \cite{AB09}).  To see this, we use the following example. Let
$\psi:[-\rho, \rho] \to \R^2$ be a convex curve given by $\psi(t)=(t, \frac{\lambda}{2}t^2)$. The
chord joining $\psi(-\rho)$ to $\psi(\rho)$ has length $2\rho$. The length-excess of
Alexander-Bishop is given by
$$\align
\int_{-\rho}^{\rho} \sqrt{1+\lambda^2t^2}\,dt -2\rho
& = \int_{-\rho}^{\rho} [\sqrt{1+\lambda^2t^2} -1]\,dt \\
& = \int_{-\rho}^{\rho} \frac{\lambda^2 t^2}{\sqrt{1+\lambda^2t^2} +1}\,dt \sim \frac{2}{3}\lambda^2\rho^3.
\endalign
$$
However, using our angular excess, we have
$$\theta(\rho)\sim \lambda\tan\rho \sim \lambda\rho$$
for sufficiently small $\rho$.

\bigskip
\noindent {\bf Acknowledgement}: Authors are grateful to Professor
Stephanie Alexander (University of Illinois at Urbana-Champaign),
Professor Vitali Kapovitch (Toronto) and Xueping Li (Capital
Normal University, Beijing) for their criticisms on an earlier
version of this paper.

\enddocument